\documentclass[leqno,11pt]{article}
\usepackage[spanish]{babel}
\usepackage[utf8]{inputenc}
\usepackage{amscd}
\usepackage{amsthm}
\usepackage{amssymb}
\usepackage{amsmath}
\usepackage{graphics}
\usepackage{graphicx}
\usepackage{verbatim}
\usepackage{color}
\newtheorem{myexa}{Example}

\newtheorem{mytheo}{Theorem}
\newtheorem{mylemma}{Lemma}

\newtheorem{myprop}{Proposition}
\newtheorem*{myrem}{{\em Remark}}

\newcommand{\pt}{\mbox{$\succ$\hspace{-1ex}$\longrightarrow$}}

\decimalpoint

\begin{document}

\begin{center}
	{\sc Optimal Bayesian Estimation of a Regression Curve, a Conditional Density and a Conditional Distribution}\vspace{2ex}\\
	A.G. Nogales\vspace{2ex}\\
	Dpto. de Matem\'aticas, Universidad de Extremadura\\
	Avda. de Elvas, s/n, 06006--Badajoz, SPAIN.\\
	e-mail: nogales@unex.es
\end{center}
\vspace{.4cm}
\begin{quote}
	\hspace{\parindent} {\small {\sc Abstract.}  		 In this paper several related estimation problems are addressed from a Bayesian point of view and optimal estimators are obtained for each of them when some natural loss functions are considered. Namely, we are interested in estimating a regression curve. Simultaneously, the estimation problems of a conditional distribution function, or a conditional density, or even the conditional distribution itself, are considered. All these problems are posed  in a sufficiently general framework to cover continuous and discrete, univariate and multivariate, parametric and non-parametric cases, without the need to use a specific prior distribution. The loss functions considered come naturally from the quadratic error loss function comonly used in estimating a real function of the unknown parameter. The cornerstone of the mentioned Bayes estimators is the posterior predictive distribution. Some examples are provided to illustrate these results. 
	}
\end{quote}

\vfill
\begin{itemize}
	\item[] \hspace*{-1cm} {\em AMS Subject Class.} (2010): {\em Primary\/} 62F15, 62G07
	{\em Secondary\/} 62Jxx
	\item[] \hspace*{-1cm} {\em Key words and phrases: } Bayesian estimation of a regression curve, posterior predictive distribution.
\end{itemize}

\newpage

\section{Introduction.}

In Statistics, the expression \emph{the probability of an event} $A$  (written $P_\theta(A)$) is, in general, ambiguous as it depends of the unknown parameter $\theta$. Before conducting the experiment, a Bayesian statistician, provided with the prior distribution, possesses  a natural candidate: the prior predictive probability of $A$ since it is the prior mean of the probabilities of $A$. However, according to Bayesian philosophy, after performing the experiment and having observed the data $\omega$,  a reasonable estimation is the posterior predicitive probability of $A$ given $\omega$ because it is the posterior mean of the probabilities of $A$ given $\omega$. It can be shown not only that this is the Bayes estimator of the probability $P_\theta(A)$ of $A$ for the squared error loss function but also that the posterior predictive distribution is the Bayes estimator of the sampling probability distribution $P_\theta$ for the squared variation total loss function and that the posterior predictive density is the Bayes estimator of its density for the $L^1$-squared loss function. Notice that these loss functions should be considered natural in the sense that they are derived directly from the quadratic error loss function commonly used in the estimation of a real function of the parameter. 
Nogales (2021) contains accurate statements and proofs of these results which are nothing but a functional generalization of Theorem 1.1 (more specifically of its Corollary 1.2.(a)) of Lehmann et al. (1998, p. 228) which obtains the Bayes estimator  of a real function of the parameter for the squared error loss function.

This article is about estimating a regression curve and some other related problems like the estimation of a conditional density or a conditional distribution function or even of the conditional distribution itself, ever from a Bayesian point of view. 
 Therefore, this paper should be considered as the conditional counterpart of Nogales (2021) and the results below become the functional extension of Lehmann et al. (1998, Theorem 1.1) for the conditional case. So it is not surprising that the posterior predictive distribution is the keystone in the estimation problems to be discussed below. Some examples are presented in Section \ref{secex} to illustrate the results. See Nogales (2021), and the references therein, for other examples of determination of the posterior predictive distribution.  But in practice, the explicit evaluation of the posterior predictive distribution could be cumbersome and its simulation may become preferable.  Gelman et al. (2014) is also a good reference for such simulation methods and, hence, for the computation of the Bayes estimators of the conditional density and the regression curve.

The posterior predictive distribution has been presented as the base for Predictive Inference, which seeks to make inferences about a new unknown observation from the previous random sample (in contrast to the greater emphasis that statistical inference makes on the estimation and contrast of parameters since its mathematical foundations  in the early twentieth century, see Geisser (1993) or Gelman et al. (2014)). 
With that idea in mind, it has also been used in other areas such as model selection, testing for discordancy, goodness of fit, perturbation analysis or classification (see addtional fields of application in Geisser (1993) and Rubin (1984)), but never as a possible solution for the Bayesian problems of estimating an unconditional or conditional density. The reader is referred to the references within Nogales (2021) for other uses of the posterior predictive distribution in Bayesian statistics.


In what follows we place ourselves in a general framewok for Bayesian inference, as described in Barra (1971).

\section{The framework.}

Let $(\Omega,\mathcal A,\{P_\theta\colon\theta\in(\Theta,\mathcal T,Q)\})$ be a Bayesian statistical experiment and $X_i:(\Omega,\mathcal A,\{P_\theta\colon\theta\in(\Theta,\mathcal T,Q)\})\rightarrow(\Omega_i,\mathcal A_i)$, $i=1,2$, two statistics. Consider the Bayesian experiment image of $(X_1,X_2)$
$$(\Omega_1\times\Omega_2,\mathcal A_1\times\mathcal  A_2,\{P_\theta^{(X_1,X_2)}\colon\theta\in(\Theta,\mathcal T,Q)\}).
$$
In the next, we shall suppose that $P^{(X_1,X_2)}(\theta,A_{12}):=P_\theta^{(X_1,X_2)}(A_{12})$, $\theta\in\Theta, \  A_{12}\in \mathcal A_1\times\mathcal  A_2$, is a Markov kernel. 
Let us write $R_\theta=P_\theta^{(X_1,X_2)}$.

The Bayesian experiment corresponding to a sample of size $n$ of the joint distribution of $(X_1,X_2)$ is 
$$\big((\Omega_1\times\Omega_2)^n,(\mathcal A_1\times\mathcal  A_2)^n,\big\{R_\theta^n\colon\theta\in(\Theta,\mathcal T,Q)\big\}\big)\qquad \hbox{(1)}
$$
We write $R^n(\theta,A'_{12,n})=R_\theta^n(A'_{12,n})$ for $A'_{12,n}\in (\mathcal A_1\times\mathcal  A_2)^n$ and 
$$\Pi_{12,n}:=Q\otimes R^n,
$$
for the joint distribution of the parameter and the sample, i.e.
$$\Pi_{12,n}(A'_{12,n}\times T)=\int_TR_\theta^n(A'_{12,n})dQ(\theta),\quad A'_{12,n}\in(\mathcal A_1\times\mathcal  A_2)^n, T\in\mathcal T.
$$
The corresponding prior predictive distribution $\beta_{12,n,Q}^*$ is
$$\beta_{12,n,Q}^*(A'_{12,n})=\int_\Theta R_\theta^n(A'_{12,n})dQ(\theta),\quad A'_{12,n}\in(\mathcal A_1\times\mathcal  A_2)^n.
$$
The posterior distribution is a Markov kernel $$R_n^*:((\Omega_1\times\Omega_2)^n,(\mathcal A_1\times\mathcal A_2)^n)\pt (\Theta,\mathcal T)
$$
such that, for all $A'_{12,n}\in(\mathcal A_1\times\mathcal  A_2)^n$ and $T\in\mathcal T$,
$$\Pi_{12,n}(A'_{12,n}\times T)=\int_TR_\theta^n(A'_{12,n})dQ(\theta)
=\int_{A'_{12,n}}R_n^*(x',T)d\beta_{12,n,Q}^*(x').
$$
Let us write $R_{n,x'}^*(T):=R_n^*(x',T)$.

The posterior predictive distribution on $\mathcal A_1\times\mathcal A_2$  is the Markov kernel
$${R_{n}^*}^{R}:((\Omega_1\times\Omega_2)^n,(\mathcal A_1\times\mathcal  A_2)^n)\pt(\Omega_1\times\Omega_2,\mathcal A_1\times\mathcal  A_2)
$$ 
defined, for $x'\in(\Omega_1\times\Omega_2)^n$, by
$${R_{n}^*}^{R}(x',A_{12}):=
\int_\Theta R_\theta(A_{12})dR_{n,x'}^*(\theta)
$$
It follows that, with obvious notations, 
$$\int_{\Omega_1\times\Omega_2}f(x)d{R_{n,x'}^*}^{\!\!\!\!R}(x)=\int_\Theta\int_{\Omega_1\times\Omega_2}f(x)dR_\theta(x)dR_{n,x'}^*(\theta)
$$
for any non-negative or integrable real random variable (r.r.v. for short) $f$. 

We can also consider the posterior predictive distribution on $(\mathcal A_1\times\mathcal A_2)^n$  defined as the Markov kernel
$${R_{n}^*}^{R^n}:((\Omega_1\times\Omega_2)^n,(\mathcal A_1\times\mathcal  A_2)^n)\pt((\Omega_1\times\Omega_2)^n,(\mathcal A_1\times\mathcal  A_2)^n)
$$ 
such that 
$${R_{n}^*}^{R^n}(x',A'_{12,n}):=
\int_\Theta R_\theta^n(A'_{12,n})dR_{n,x'}^*(\theta)
$$

According to Theorem 1 of Nogales (2021), it is the Bayes estimator of the distribution $R_\theta^n$ for the squared total variation function, i.e.,
\begin{gather*}
	\int_{(\Omega_1\times\Omega_2)^n\times\Theta}\sup_{A_{12,n}\in(\mathcal A_1\times\mathcal A_2)^n}\left|{R_{n,x'}^*}^{\!\!\!\!R^n}(A_{12,n})-R_\theta^n(A_{12,n})\right|^2d\Pi_{12,n}(x',\theta)\le\\
	\int_{(\Omega_1\times\Omega_2)^n\times\Theta}\sup_{A_{12,n}\in(\mathcal A_1\times\mathcal A_2)^n}\left| M(x',A_{12,n})-R_\theta^n(A_{12,n})\right|^2d\Pi_{12,n}(x',\theta),
\end{gather*}
for every Markov kernel $M:(\Omega_1\times\Omega_2,\mathcal A_1\times\mathcal A_2)^n\pt (\Omega_1\times\Omega_2,\mathcal A_1\times\mathcal A_2)^n$.

It can be readily checked that 
$$\left[{R_{n,x'}^*}^{\!\!\!\!R^n}\right]^{\pi'_1}={R_{n,x'}^*}^{\!\!\!\!R},
$$
where $\pi'_1(x'):=x'_1:=(x'_{11},x'_{21})$ for $x'\in(\Omega_1\times\Omega_2)^n$. So, Theorem 2 of Nogales (2021) shows that 
\begin{gather*}
	\int_{(\Omega_1\times\Omega_2)^n\times\Theta}\sup_{A_{12}\in\mathcal A_1\times\mathcal A_2}\left|{R_{n,x'}^*}^{\!\!\!\!R}(A_{12})-R_\theta(A_{12})\right|^2d\Pi_{12,n}(x',\theta)\le\\
	\int_{(\Omega_1\times\Omega_2)^n\times\Theta}\sup_{A_{12}\in\mathcal A_1\times\mathcal A_2}\left| M(x',A_{12})-R_\theta(A_{12})\right|^2d\Pi_{12,n}(x',\theta),
\end{gather*}
for every Markov kernel $M:(\Omega_1\times\Omega_2,\mathcal A_1\times\mathcal A_2)^n\pt (\Omega_1\times\Omega_2,\mathcal A_1\times\mathcal A_2)$.

We introduce some notations for $(x',x,\theta)\in(\Omega_1\times\Omega_2)^n\times(\Omega_1\times\Omega_2)\times\Theta$:
\begin{gather*}
\pi'(x',x,\theta):=x',\quad  \pi'_i(x',x,\theta):=x'_i:=(x'_{i1},x'_{i2}),\quad 1\le i\le n,\\
\pi(x',x,\theta):=x,\quad \pi_i(x',x,\theta):=x_i,\quad i=1,2,\\
q(x',x,\theta):=\theta.
\end{gather*}
Let us consider the probability space  
\begin{gather*}
((\Omega_1\times\Omega_2)^n\times(\Omega_1\times\Omega_2)\times\Theta,(\mathcal A_1\times\mathcal A_2)^n\times(\mathcal A_1\times\mathcal A_2)\times\mathcal T,\Pi_n),\qquad \hbox{(2)}
\end{gather*}
where
$$\Pi_n(A'_{12,n}\times A_{12}\times T)=\int_T R_\theta(A_{12}) R_\theta^n(A'_{12,n})dQ(\theta),
$$
when $A'_{12,n}\in(\mathcal A_1\times\mathcal A_2)^n$, $A_{12}\in \mathcal A_1\times\mathcal A_2$ and $T\in\mathcal T$. 

So, for a r.r.v. $f$ on $((\Omega_1\times\Omega_2)^n\times(\Omega_1\times\Omega_2)\times\Theta,(\mathcal A_1\times\mathcal A_2)^n\times(\mathcal A_1\times\mathcal A_2)\times\mathcal T)$,
$$\int f d\Pi_n=\int_\Theta\int_{(\Omega_1\times\Omega_2)^{n}}
\int_{(\Omega_1\times\Omega_2)} f(x',x,\theta)dR_\theta(x)dR_\theta^n(x')dQ(\theta)\qquad {\rm (3)}
$$
provided that the integral exists. 
Moreover, for a r.r.v. $h$ on $((\Omega_1\times\Omega_2)\times\Theta,(\mathcal A_1\times\mathcal A_2)\times\mathcal T)$,
$$\int h d\Pi_n=\int_\Theta\int_{\Omega_1\times\Omega_2}h(x,\theta)dR_\theta(x)dQ(\theta)
=\int_{\Omega_1\times\Omega_2}\int_\Theta h(x,\theta)dR^*_{1,x}(\theta)d\beta^*_{12,1,Q}(x).
$$ 

The following proposition is straightforward.

\begin{myprop}\begin{gather*}
\Pi_n^{(\pi',q)}(A'_{12.n}\times T)=\Pi_{12,n}(A'_{12,n}\times T)=
\int_TR_\theta^n(A'_{12,n})dQ(\theta)=
\int_{A'_{12,n}} R^{*}_{n,x'}(T)d\beta_{12,n,Q}^*(x')\\
\Pi_n^{(\pi,q)}(A_{12}\times T)=\Pi_{12,1}(A_{12}\times T)=
\int_TR_\theta(A_{12})dQ(\theta)=
\int_{A_{12}} R^{*}_{1,x}(T)d\beta_{12,1,Q}^*(x)\\
\Pi_n^q=Q,\quad \Pi_n^{(\pi',q)}=\Pi_{12,n},\quad \Pi_n^{\pi'}=\beta_{12,n,Q}^*,\quad 
\Pi_n^{(\pi,q)}=\Pi_{12,1},\quad
\Pi_n^{\pi}=\beta_{12,1,Q}^*\\ \Pi_n^{\pi'|q=\theta}=R_\theta^n,\quad \Pi_n^{\pi|q=\theta}=R_\theta,\quad \Pi_n^{q|\pi'=x'}=R^*_{n,x'},\quad \Pi_n^{q|\pi=x}=R^*_{1,x'}\\
P_\theta^{X_1}=R_\theta^{\pi_1},\quad P_\theta^{X_2|X_1=x_1}=R_\theta^{\pi_2|\pi_1=x_1},\\
E_{P_\theta}(X_2|X_1=x_1)=E_{R_\theta}(\pi_2|\pi_1=x_1).
\end{gather*}\end{myprop}

In particular, the probability space (2) contains all the basic ingredients of the Bayesian experiment (1), i.e., the prior distribution, the sampling probabilities, the posterior distributions and the prior predictive distribution. In addition, it becomes the natural framework to address the estimation problems of this paper, as we can see in the next.

\section{Bayes estimator of the conditional distribution.}

An estimator of the conditional distribution $P_\theta^{X_2|X_1}$ from a $n$-sized sample of the joint distribution of $(X_1,X_2)$ is a Markov kernel 
$$M:((\Omega_1\times\Omega_2)^n\times\Omega_1,(\mathcal A_1\times\mathcal A_2)^n\times\mathcal A_1)\pt (\Omega_2,\mathcal A_2)
$$
so that, being observed $x'=((x'_{11},x'_{21}),\dots,(x'_{1n},x'_{2n}))\in (\Omega_1\times \Omega_2)^n$, $M(x',x_1,\cdot)$ is a probability measure on $\mathcal A_2$ that is considered as an estimation of the conditional distribution $P_\theta^{X_2|X_1=x_1}$ for a given $x_1\in\Omega_1$.

From a Bayesian point of view, the Bayes estimator of the conditional distribution $P^{X_2|X_1}=R^{\pi_2|\pi_1}$
 is a Markov kernel
$$M:((\Omega_1\times\Omega_2)^n\times\Omega_1,(\mathcal A_1\times\mathcal A_2)^n\times\mathcal A_1)\pt (\Omega_2,\mathcal A_2)
$$
minimizing the Bayes risk
\begin{gather*}
\int_{(\Omega_1\times\Omega_2)^n\times\Theta}\int_{\Omega_1}
\sup_{A_2\in\mathcal A_2}\big|M(x',x_1,A_2)-R_\theta^{\pi_2|\pi_1=x_1}(A_2)\big|^2
dR_\theta^{\pi_1}(x_1)d\Pi_{12,n}(x',\theta)=\\
\int_\Theta\int_{(\Omega_1\times\Omega_2)^n}\int_{\Omega_1}
\sup_{A_2\in\mathcal A_2}\big|M(x',x_1,A_2)-R_\theta^{\pi_2|\pi_1=x_1}(A_2)\big|^2
dR_\theta^{\pi_1}(x_1)dR_\theta^n(x')dQ(\theta)=\\
\int_{(\Omega_1\times\Omega_2)^n\times(\Omega_1\times\Omega_2)\times\Theta}
\sup_{A_2\in\mathcal A_2}\big|M(x',x_1,A_2)-R_\theta^{\pi_2|\pi_1=x_1}(A_2)\big|^2
d\Pi_n(x',x,\theta).
\end{gather*}

The next result yields the Bayes estimator of the conditional distribution $P_\theta^{X_2|X_1}$ from the posterior predictive distribution.

\begin{mytheo}\label{Theo1}\rm Suppose that the $\sigma$-field $\mathcal A_2$ is separable. Then, the conditional distribution of $\pi_2$ given $\pi_1=x_1$ with respect to the posterior predictive distribution ${R^*_{n,x'}}^{\!\!\!\!R}$, 
	$$M_n^*(x',x_1,A_2):=\left({R^*_{n,x'}}^{\!\!\!\!R}\right)^{\pi_2|\pi_1=x_1}(A_2),
	$$
	is the Bayes estimator of the conditional distribution   $R^{\pi_2|\pi_1}$ for the squared total variation loss function, i.e.,
\begin{gather*}
\int_{(\Omega_1\times\Omega_2)^{n+1}\times\Theta}
\sup_{A_2\in\mathcal A_2}\big|M_n^*(x',x_1,A_2)-R_\theta^{\pi_2|\pi_1=x_1}(A_2)\big|^2
d\Pi_n(x',x,\theta)\le
\\
\int_{(\Omega_1\times\Omega_2)^{n+1}\times\Theta}
\sup_{A_2\in\mathcal A_2}\big|M(x',x_1,A_2)-R_\theta^{\pi_2|\pi_1=x_1}(A_2)\big|^2
d\Pi_n(x',x,\theta)
\end{gather*}
for any estimator $M$ of the conditional distribution $R^{\pi_2|\pi_1}$.
\end{mytheo}

	Fix an event $A_2\in\mathcal A_2$ and define
	$$H_{A_2}(x',x_1,\theta):=P_\theta^{X_2|X_1=x_1}(A_2)=R_\theta^{\pi_2|\pi_1=x_1}(A_2).
	$$
	Jensen's inequality could help to achieve a proof of the theorem  if we can prove the following result.
	\begin{mylemma}	\label{lem1}
	$$E_{\Pi_n}(H_{A_2}|\pi'=x',\pi_1=x_1)=M_n^*(x',x_1,A_2),
	$$
	\end{mylemma}

	\section{Bayes estimator of the conditional density.}

When the joint distribution $R_\theta=P_\theta^{(X_1,X_2)}$ has a density $f_\theta$ with respect to the product of two $\sigma$-finite measures $\mu_1$ and $\mu_2$ on $\mathcal A_1$ and $\mathcal A_2$, resp., the conditional density is
$$f_\theta^{X_2|X_1=x_1}(x_2):=\frac{f_\theta(x_1,x_2)}{f_{\theta,X_1}(x_1)}
$$
for almost every $x_1$, where $f_{\theta,X_1}(x_1)$ stands for the marginal density of $X_1$. 

An estimator of the conditional density $f_\theta^{X_2|X_1}$ from a $n$-sized sample of the joint distribution of $(X_1,X_2)$ is a map 
$$m:((\Omega_1\times\Omega_2)^n\times\Omega_1\times\Omega_2,(\mathcal A_1\times\mathcal A_2)^n\times\mathcal A_1\times\mathcal A_2)\rightarrow (\mathbb R,\mathcal R)
$$
so that, being observed $x'=((x'_{11},x'_{21}),\dots,(x'_{1n},x'_{2n}))\in (\Omega_1\times \Omega_2)^n$, $m(x',x_1,\cdot)$  is considered as an estimation of the conditional density $f_\theta^{X_2|X_1=x_1}$ of $X_2$ given $X_1=x_1$.

It is well known (see, for instance, Ghosal et al. (2017, p.126)) that, given two probability measures $P_1$ and $P_2$ on a measurable space $(\Omega,\mathcal A)$ having densities $p_1$ and $p_2$ with respect to a $\sigma$-finite measure $\mu$, 
$$\sup_{A\in\mathcal A}|P_1(A)-P_2(A)|=\frac12\int_\Omega|p_1-p_2|d\mu.
$$
So the Bayesian estimation of the conditional distribution $P_\theta^{X_2|X_1=x_1}=R_\theta^{\pi_2|\pi_1=x_1}$ for the squared total variation loss function corresponds to the Bayesian estimation of its density $f_\theta^{X_2|X_1=x_1}$ for the $L^1$-squared loss function. 
Hence, according to Theorem \ref{Theo1}, the Bayes estimator of the conditional density $f_\theta^{X_2|X_1=x_1}$ for the $L^1$-squared loss function is the $\mu_2$-density ${f^*_{n,x'}}^{\!\!\!\!X_2|X_1=x_1}$ of the conditional distribution
$$\left[{R^{*}_{n,x'}}^{\!\!\!\!R}\right]^{\pi_2|\pi_1=x_1}.
$$
Notice that 
\begin{gather*}
	{R^{*}_{n,x'}}^{\!\!\!\!R}(A_1\times A_2)=\int_\Theta R_\theta(A_1\times A_2)dR^*_{n,x'}(\theta)=\\
	\int_\Theta\int_{A_1\times A_2}f_\theta(x_1,x_2)d(\mu_1\times\mu_2)(x_1,x_2)r^*_{n,x'}(\theta)dQ(\theta)=\\
	\int_{A_1\times A_2}\int_\Theta f_\theta(x_1,x_2) r^*_{n,x'}(\theta)dQ(\theta)d(\mu_1\times\mu_2)(x_1,x_2)
\end{gather*}
where $r^*_{n,x'}(\theta)$ denotes the $Q$-density of the posterior distribution $R^*_{n,x'}$. So $r^*_{n,x'}(\theta)$ is of the form $K(x')f_{n,\theta}(x')$, where
$$f_{n,\theta}(x'):=\prod_{i=1}^nf_\theta(x'_i)
$$
is the $(\mu_1\times\mu_2)^n$-density of $R_\theta^n$. Hence the $\mu_1\times\mu_2$-density of the posterior predictive distribution ${R^{*}_{x'}}^{\!\!R}$ is
$$f^*_{n,x'}(x_1,x_2):=\int_\Theta f_\theta(x_1,x_2) r^*_{n,x'}(\theta)dQ(\theta),
$$
and its first marginal is
$$f^*_{n,x',1}(x_1):=\int_{\Omega_2}\int_\Theta f_\theta(x_1,t) r^*_{n,x'}(\theta)dQ(\theta)d\mu_2(t).
$$
So, we have proved the following result. 

\begin{mytheo}\rm Suppose that $\mathcal A_2$ is separable. The Bayes estimator of the conditional density $f_\theta^{X_2|X_1}$ for the $L^1$-squared loss function is the  $\mu_2$-density $${f^*_{n,x'}}^{\!\!\!\!X_2|X_1=x_1}(x_2):=\frac{f^*_{n,x'}(x_1,x_2)}{f^*_{n,x',1}(x_1)}=\frac{\int_\Theta f_\theta(x_1,x_2) r^*_{n,x'}(\theta)dQ(\theta)}{\int_{\Omega_2}\int_\Theta f_\theta(x_1,t) r^*_{n,x'}(\theta)dQ(\theta)d\mu_2(t)}
$$ 
of the conditional distribution
	$\left[{R^{*}_{n,x'}}^{\!\!\!\!R}\right]^{\pi_2|\pi_1}$ of $\pi_2$ given $\pi_1$ with respect to the posterior predictive distribution ${R^{*}_{n,x'}}^{\!\!\!\!R}$,
	i.e.
\begin{gather*}
	\int_{(\Omega_1\times\Omega_2)^{n+1}\times\Theta}
	\left(\int_{\Omega_2}\left|{f^*_{n,x'}}^{\!\!\!\!X_2|X_1=x_1}(t)-f_\theta^{X_2|X_1=x_1}(t)\right|d\mu_2(t)\right)^2d\Pi_n(x',x,\theta)\le\\
	\int_{(\Omega_1\times\Omega_2)^{n+1}\times\Theta}
	\left(\int_{\Omega_2}\left|m(x',x_1,t)-f_\theta^{X_2|X_1=x_1}(t)\right|d\mu_2(t)\right)^2d\Pi_n(x',x,\theta),
\end{gather*}
for any estimator $m$ of the conditional density.
	\end{mytheo}

	\section{Bayes estimator of the conditional distribution function.}

When $X_2$ is a r.r.v. we may be interested in the estimation of the conditional distribution function of $X_2$ given $X_1=x_1$
$$F_\theta(x_1,t):=P_\theta(X_2\le t|X_1=x_1)=R_\theta^{\pi_2|\pi_1=x_1}(]-\infty,t]).
$$
An estimator of such a conditional distribution function from a $n$-sized sample of $R_\theta$ is a map of the form
$$F:(x',x_1,t)\in (\Omega_1\times \mathbb R)^n\times\Omega_1\times\mathbb R\mapsto 
F(x',x_1,t):=M(x',x_1,]-\infty,t])\in[0,1]
$$
for a Markov kernel 
$$M:
((\Omega_1\times \mathbb R)^n\times\Omega_1,(\mathcal A_1\times\mathcal R)^n \times\mathcal A_1)
\pt 
(\mathbb R,\mathcal R)
$$
An optimal estimator of the conditional distribution function $F_\theta$ for the $L^\infty$-squared loss function from a Bayesian point of view (i.e. a Bayes estimator) is an estimator $F_n^*$ minimizing the Bayes risk, i.e., such that
\begin{gather*}
	\int_{(\Omega_1\times\mathbb R)^{n+1}\times\Theta}\sup_{t\in\mathbb R}|F_n^*(x',x_1,t)-F_\theta(x_1,t)|^2d{\Pi_n}(x',x,\theta)
	\le\\
	\int_{(\Omega_1\times\mathbb R)^{n+1}\times\Theta}\sup_{t\in\mathbb R}|F(x',x_1,t)-F_\theta(x_1,t)|^2d{\Pi_n}(x',x,\theta)
\end{gather*}
for any estimator $F$ of the conditional distribution function $F_\theta$.

A natural candidate is the conditional distribution function for the posterior predictive distribution as it is stated in the next theorem.

\begin{mytheo}\rm \label{Theo3} The posterior predictive conditional distribution function 
	$$F_n^*(x',x_1,t):= \left({R^*_{n,x'}}^{\!\!\!\!R}\right)^{\pi_2|\pi_1=x_1}\!\!(]-\infty,t])
	$$ 
	is the Bayes estimator of the conditional distribution function $F_\theta$ for the $L^\infty$-squared loss function.
\end{mytheo}

	\section{Bayes estimator of a regression curve.}

Now suppose that $X_2$ is an squared-integrable r.r.v. So $(\Omega_2,\mathcal A_2)=(\mathbb R,\mathcal R)$. The regression curve of $X_2$ given $X_1$ is the map $x_1\in\Omega_1\mapsto r_\theta(x_1):=E_\theta(X_2|X_1=x_1)$. An estimator of the regression curve $r_\theta$ from a sample of size $n$ of the joint distribution of $(X_1,X_2)$ is a statistic
$$m:(x',x_1)\in(\Omega_1\times \mathbb R)^n\times\Omega_1\longmapsto m(x',x_1)\in\mathbb R,
$$
so that, being observed $x'\in (\Omega_1\times \mathbb R)^n$, $m(x',\cdot)$ is the estimation of $r_\theta$.

From a frequentist point of view, the simplest way to evaluate the error in estimating an unknown regression curve is to use the expectation of the quadratic deviation (see Nadaraya (1989, p. 120)):
\begin{gather*}
	E_\theta\left[\int_{\Omega_1}(m(x',x_1)-r_\theta(x_1))^2dP_\theta^{X_1}(x_1)\right]=\\
	\int_{(\Omega_1\times \mathbb R)^n}\int_{\Omega_1}(m(x',x_1)-r_\theta(x_1))^2dR_\theta^{\pi_1}(x_1) dR_\theta^n(x').
\end{gather*}

From a Bayesian point of view, the Bayes estimator of the regression curve $r_\theta$ should minimize the Bayes risk (i.e. the prior mean of the expectation of the quadratic deviation)
\begin{gather*}
	\int_\Theta\int_{(\Omega_1\times \mathbb R)^n}\int_{\Omega_1}(m(x',x_1)-r_\theta(x_1))^2dR_\theta^{\pi_1}(x_1) dR_\theta^n(x')dQ(\theta)=\\
	E_{\Pi_n}\big[(m(x',x_1)-r_\theta(x_1))^2\big].
\end{gather*}

The following result solves the problem of estimation of the regression curve from a Bayesian point of view.

\begin{mytheo}\rm \label{Theo4} The regression curve of $\pi_2$ on $\pi_1$ with respect to the posterior predictive distribution ${R^*_{n,x'}}^{\!\!\!\!\!R}$ $$m_n^*(x',x_1):=E_{{R^*_{n,x'}}^{\!\!\!\!\!R}}(\pi_2|\pi_1=x_1)$$ is the Bayes estimator of the regression curve $r_\theta(x_1):=E_\theta(X_2|X_1=x_1)$ for the squared error loss function, i.e., $$E_{\Pi_n}[(m_n^*(x',x_1)-r_\theta(x_1))^2]\le 
	E_{\Pi_n}[(m_n(x',x_1)-r_\theta(x_1))^2]
	$$
	for any other estimator $m_n$ of the regression curve $r_\theta$.	
	\end{mytheo}

\begin{myrem}\rm (Estimation of the regression curve when densities are available) According to the previous results, when $R_\theta$ has density $f_\theta$ with respect to the product $\mu_1\times\mu_2$ of two $\sigma$-finite measures, the $\mu_2$-density ${f^*_{n,x'}}^{\!\!\!\!X_2|X_1=x_1}$  of the conditional distribution $\left[{R^{*}_{n,x'}}^{\!\!\!\!R}\right]^{\pi_2|\pi_1=x_1}$ is	
$${f^*_{n,x'}}^{\!\!\!\!X_2|X_1=x_1}(x_2):=\frac{f^*_{n,x'}(x_1,x_2)}{f^*_{n,x',1}(x_1)}=\frac{\int_\Theta f_\theta(x_1,x_2) r^*_{n,x'}(\theta)dQ(\theta)}{\int_{\mathbb R}\int_\Theta f_\theta(x_1,t) r^*_{n,x'}(\theta)dQ(\theta)d\mu_2(t)},
$$
which is the Bayes estimator of the conditional density $f_\theta^{X_2|X_1}$. Hence, the Bayes estimator of the regression curve can be computed as
\begin{gather*}
m_n^*(x',x_1):=E_{{R^*_{n,x'}}^{\!\!\!\!\!R}}(\pi_2|\pi_1=x_1)=\\
\int_{\mathbb R}x_2\cdot {f^*_{n,x'}}^{\!\!\!\!X_2|X_1=x_1}(x_2)d\mu_2(x_2)=\\
\frac{\int_{\mathbb R}x_2\int_\Theta f_\theta(x_1,x_2) r^*_{n,x'}(\theta)dQ(\theta)d\mu_2(x_2)}{\int_{\mathbb R}\int_\Theta f_\theta(x_1,x_2) r^*_{n,x'}(\theta)dQ(\theta)d\mu_2(x_2)}.
\end{gather*}
	\end{myrem}

\section{Examples.}\label{secex}

\begin{myexa}\rm Let us suppose that, for $\theta,\lambda,x_1>0$, $P_\theta^{X_1}= G(1,\theta^{-1})$, $P_\theta^{X_2|X_1=x_1}=G(1,(\theta x_1)^{-1})$ and $Q=G(1,\lambda^{-1})$, where $G(\alpha,\beta)$ denotes the gamma distribution of parameters $\alpha,\beta>0$. Hence the joint density of $X_1$ and $X_2$ is
	$$f_\theta(x_1,x_2)=\theta^2 x_1\exp\{-\theta x_1(1+x_2)\}I_{]0,\infty[^2}(x_1,x_2).
	$$
	So, the density of $R_\theta^n$ is
	$$f_{n,\theta}(x')=\theta^{2n}\cdot \prod_{i=1}^nx'_{i1}\cdot\exp\big\{-\theta\textstyle\sum_{i=1}^nx'_{i1}(1+x'_{i2})\big\}\cdot I_{]0,\infty[^{2n}}(x'),
	$$
	and the posterior $Q$-density given $x'$ is
	$$\frac{dR^*_{n,x'}(\theta)}{dQ}=:r^*_{n,x'}(\theta)=K(x')f_{n,\theta}(x')
	$$
	where $K(x')=[\int_0^\infty f_{n,\theta}(x')dQ(\theta)]^{-1}$. 
	
	So, the posterior predictive density given $x'\in ]0,\infty[^{2n}$ is
	\begin{gather*}
	f^*_{n,x'}(x)=\int_0^\infty f_\theta(x)r^*_{n,x'}(\theta)dQ(\theta)=\\
    \lambda K(x') x_1\prod_{i=1}^nx'_{i1}\int_0^\infty \theta^{2n+2} \exp\left\{-\theta[\lambda+x_1(1+x_2)+\textstyle\sum_{i=1}^nx'_{i1}(1+x'_{i2})]\right\}d\theta\cdot I_{]0,\infty[^2}(x).
	\end{gather*}
Since
$$\int_0^\infty \theta^n\exp\{-a\theta\}d\theta=\frac{n!}{a^{n+1}},
$$
we have that
$$f^*_{n,x'}(x)=\frac{(2n+2)!\lambda K(x') x_1\prod_{i=1}^nx'_{i1}}{[\lambda+x_1(1+x_2)+\textstyle\sum_{i=1}^nx'_{i1}(1+x'_{i2})]^{(2n+3)}}
$$
and its first marginal is
$$f^*_{n,x',1}(x_1)=\int_{\mathbb R} f^*_{n,x'}(x_1,x_2)dx_2=\int_0^\infty\frac{A}{(Bt+C)^m}dt=\frac{A}{(m-1)BC^{m-1}}
$$
where
$$m=2n+3,\quad A=(2n+2)!\lambda K(x') x_1\prod_{i=1}^nx'_{i1},\quad B=x_1,\quad\hbox{and}\quad C=\lambda+x_1+\sum_{i=1}^nx'_{i1}(1+x'_{i2}).
$$
So
$$f^*_{n,x',1}(x_1)=\frac{(2n+2)!\lambda K(x') x_1\prod_{i=1}^nx'_{i1}}{(2n+2)x_1[\lambda+x_1+\sum_{i=1}^nx'_{i1}(1+x'_{i2})]^{2n+2}}.
$$
The Bayes estimator of the conditional density $f_\theta^{X_2|X_1=x_1}(x_2)=\theta x_1\exp\{-\theta x_1x_2\}I_{]0,\infty[}(x_2)$ is, for $x_1,x_2>0$, 
\begin{gather*}
{f^*_{n,x'}}^{\!\!\!\!X_2|X_1=x_1}(x_2):=\frac{f^*_{n,x'}(x_1,x_2)}{f^*_{n,x',1}(x_1)}=\\(2n+2)x_1\frac{[\lambda+x_1+\sum_{i=1}^nx'_{i1}(1+x'_{i2})]^{2n+2}}{[\lambda+x_1(1+x_2)+\sum_{i=1}^nx'_{i1}(1+x'_{i2})]^{2n+3}}=\\
\frac{(2n+2)x_1a_n(x',x_1)^{2n+2}}{(x_1x_2+a_n(x',x_1))^{2n+3}}
\end{gather*}
where 
$$a_n(x',x_1)=\lambda+x_1+\sum_{i=1}^nx'_{i1}(1+x'_{i2}).
$$
The Bayes estimator of the conditional distribution function
$$F_\theta(x_1,t):=P_\theta(X_2\le t|X_1=x_1)
$$
is, for $t>0$,
\begin{gather*}\begin{split}
F_n^*(x',x_1,t)&=\int_0^t
(2n+2)x_1\frac{[\lambda+x_1+\sum_{i=1}^nx'_{i1}(1+x'_{i2})]^{2n+2}}{[\lambda+x_1(1+x_2)+\sum_{i=1}^nx'_{i1}(1+x'_{i2})]^{2n+3}}
dx_2\\
&=\int_0^t\frac{(2n+2)x_1a_n(x',x_1)^{2n+2}}{(x_1x_2+a_n(x',x_1))^{2n+3}}dx_2\\
&=a_n(x',x_1)^{2n+2}\left(\frac1{a_n(x',x_1)^{2n+2}}-\frac1{(x_1t+a_n(x',x_1))^{2n+2}}\right)
\\&=1-\left(1+\frac{x_1t}{a_n(x',x_1)}\right)^{-2n-2}.
\end{split}\end{gather*}

The Bayes estimator of the regression curve $r_\theta(x_1):=E_\theta(X_2|X_1=x_1)=\frac1{\theta x_1}$ is, for $x_1>0$, 
$$m^*_n(x',x_1)=\int_0^\infty x_2\cdot{f^*_{n,x'}}^{\!\!\!\!X_2|X_1=x_1}(x_2)dx_2=\frac{\lambda+x_1+\sum_{i=1}^nx'_{i1}(1+x'_{i2})}{(2n+1)x_1}=\frac{a_n(x',x_1)}{(2n+1)x_1}.
$$
	\end{myexa} 
	
\begin{myexa}\rm Let us suppose that $X_1$ has a Bernoulli distribution of unknown parameter $\theta\in]0,1[$ (i.e. $P_\theta^{X_1}=Bi(1,\theta)$) and, given $X_1=k_1\in\{0,1\}$, $X_2$ has distribution $Bi(1,1-\theta)$ when $k_1=0$ and $Bi(1,\theta)$ when $k_1=1$, i.e. $P_\theta^{X_2|X_1=k_1}=Bi(1,k_1+(1-2k_1)(1-\theta))$. We can think of tossing a coin with probability $\theta$ of getting heads ($=1$) and making a second toss of this coin if it comes up heads on the first toss, or tossing a second coin with probability $1-\theta$ of making heads if the first toss is tails ($=0$). Consider the uniform distribution on $]0,1[$ as the prior distribution $Q$.

So, the joint probability function of $X_1$ and $X_2$ is
\begin{gather*}\begin{split}
f_\theta(k_1,k_2)&=\theta^{k_1}(1-\theta)^{1-k_1}[k_1+(1-2k_1)(1-\theta)]^{k_2}[1-k_1-(1-2k_1)(1-\theta)]^{1-k_2}\\
&=\begin{cases}
\theta(1-\theta)\hbox{ if } k_2=0,\\ (1-\theta)^2\hbox{ if } k_1=0, k_2=1,\\ \theta^2\hbox{ if } k_1=1, k_2=1.
\end{cases}
\end{split}\end{gather*}
The probability function of $R_\theta^n$ is
$$f_{n,\theta}(k')=\prod_{i=1}^nf_\theta(k'_i)=\theta^{a_n(k')}(1-\theta)^{b_n(k')}
$$
for $k'=(k'_1,\dots,k'_n)=(k'_{11},k'_{12},\dots,k'_{n1},k'_{n2})\in\{0,1\}^{2n}$, where
$$a_n(k')=n_{+0}(k')+2n_{11}(k'),\quad b_n(k')=n_{+0}(k')+2n_{01}(k'),
$$
being $n_{j_1j_2}(k')$ the number of indices $i\in\{1,\dots,n\}$ such that $(k'_{i1},k'_{i2})=(j_1,j_2)$ and $n_{+j}=n_{0j}+n_{1j}$ for $j=0,1$. Notice that $a_n(k')+b_n(k')=2n$.

Hence, the posterior $Q$-density given $k'$ is
$$r^*_{n,k'}(\theta)=\frac{dR^*_{n,k'}}{dQ}(\theta)=
K(k')f_{n,\theta}(k')=K(k')\theta^{a_n(k')}(1-\theta)^{b_n(k')}
$$
where 
$$K(k')=\left[\int_0^1 f_{n,\theta}(k')dQ(\theta)\right]^{-1}=\frac1{B(a_n(k')+1,b_n(k')+1)},
$$
where $B(\alpha,\beta)$ stands for the Beta function. 

	So, the posterior predictive density given $k' \in\{0,1\}^{2n}$ is
\begin{gather*}\begin{split}
	f^*_{n,k'}(k_1,k_2)&=\int_0^1 f_\theta(k_1,k_2)r^*_{n,k'}(\theta)dQ(\theta)\\
	&=
	K(k')\int_0^1\theta^{a_n(k')+k_1}(1-\theta)^{b_n(k')+k_2}d\theta\\
	&=K(k')B(a_n(k')+k_1+1,b_n(k')+k_2+1),
	\end{split}
\end{gather*}
and its first marginal is
\begin{gather*}
		f^*_{n,k',1}(k_1)=K(k')[B(a_n(k')+k_1+1,b_n(k')+1)+B(a_n(k')+k_1+1,b_n(k')+2)].
\end{gather*}
The Bayes estimator of the conditional probability function $$f_\theta^{X_2|X_1=k_1}(k_2)=[k_1+(1-2k_1)(1-\theta)]^{k_2}[1-k_1-(1-2k_1)(1-\theta)]^{1-k_2}
$$ 
is 
\begin{gather*}
	{f^*_{n,x'}}^{\!\!\!\!X_2|X_1=k_1}(k_2):=\frac{f^*_{n,k'}(k_1,k_2)}{f^*_{n,k',1}(k_1)}=
	\begin{cases}
	\frac{2n+2}{2n+n_{+0}(k')+2n_{01}(k')+3}&\hbox{ if } k_1=k_2=0,\vspace{1ex}\\
	\frac{n_{+0}(k')+2n_{01}(k')+1}{2n+n_{+0}(k')+2n_{01}(k')+3}&\hbox{ if } k_1=0, k_2=1,\vspace{1ex}\\
	\frac{2n+3}{2n+n_{+0}(k')+2n_{01}(k')+4}&\hbox{ if } k_1=1, k_2=0,\vspace{1ex}\\
	\frac{n_{+0}(k')+2n_{01}(k')+1}{2n+n_{+0}(k')+2n_{01}(k')+4}&\hbox{ if } k_1=k_2=1.
		\end{cases}	
\end{gather*}
The Bayes estimator of the conditional mean $r_\theta(k_1):=E_\theta(X_2|X_1=k_1)=\theta^{k_1}(1-\theta)^{1-k_1}$ is, for $k_1=0,1$,
$$m^*_n(k',k_1)={f^*_{n,k'}}^{\!\!\!\!X_2|X_1=k_1}(1)=\begin{cases}
		\frac{n_{+0}(k')+2n_{01}(k')+1}{2n+n_{+0}(k')+2n_{01}(k')+3}&\hbox{ if } k_1=0,\vspace{1ex}\\
	\frac{n_{+0}(k')+2n_{01}(k')+1}{2n+n_{+0}(k')+2n_{01}(k')+4}&\hbox{ if } k_1=1.
\end{cases}	
$$
	\end{myexa}

\begin{myexa}\rm 
Let $(X_1,X_2)$ have bivariate normal distribution  with density
\begin{gather*}
f_\theta(x):=\frac1{2\pi \sigma^2\sqrt{1-\rho^2}}\exp\left\{-\frac1{2\sigma^2(1-\rho^2)}[(x_1-\theta)^2-2\rho(x_1-\theta)(x_2-\theta)+(x_2-\theta)^2]\right\}\\
=\frac1{2\pi \sigma^2\sqrt{1-\rho^2}}\exp\left\{-\frac1{2\sigma^2(1-\rho^2)}[x_1^2+x_2^2-2\rho x_1x_2-2(1-\rho)(x_1+x_2)\theta+2(1-\rho)\theta^2]
\right\},
\end{gather*}
where $\sigma>0$ and $\rho\in]-1,1[$ are supposed to be kown. 
So 
\begin{gather*}
	R_\theta=N_2\left(
	\left(\begin{array}{c}
		\theta\\
		\theta
		\end{array}\right),\sigma^2
	\left(\begin{array}{cc}
		1 & \rho\\
		\rho & 1
	\end{array}\right)
	\right),\qquad X_1,X_2\sim_\theta N(\theta,\sigma^2)\\
	P_\theta^{X_2|X_1=x_1}=N((1-\rho)\theta+\rho x_1,\sigma^2\sqrt{1-\rho^2}),\qquad
	E_\theta(X_2|X_1=x_1)=(1-\rho)\theta+\rho x_1.
	\end{gather*}
Hence, for $x'=(x'_1,\dots,x'_n)=(x'_{11},x'_{12},\dots,x'_{n1},x'_{n2})\in(\mathbb R^2)^n$,
	\begin{gather*}
f_{n,\theta}(x')=\prod_{i=1}^nf_\theta(x'_i)=\\
\frac1{[2\pi\sigma^2\sqrt{1-\rho^2}]^n}\exp\left\{
-\frac1{2\sigma^2(1-\rho^2)}\textstyle\sum_{i=1}^n\left[(x'_{i1}-\theta)^2-2\rho(x'_{i1}-\theta)(x'_{i2}-\theta)+(x'_{i2}-\theta)^2\right]\right\}.
	\end{gather*}
Let us consider the prior distribution $Q=N(\mu,\tau^2)$ whose density is 
$$g(\theta)=\frac1{\tau\sqrt{2\pi}}\exp\left\{-\frac1{2\tau^2}(\theta-\mu)^2\right\}.
$$
The posterior $Q$-density is
\begin{gather*}r^*_{n,x'}(\theta):=\frac{dR^*_{n,x'}}{dQ}(\theta)=K_1(x')f_{n,\theta}(x')=\\
\frac{K_1(x')}{[2\pi\sigma^2\sqrt{1-\rho^2}]^n}
\exp\left\{
-\frac1{2\sigma^2(1-\rho^2)}[2n(1-\rho)\theta^2-2(1-\rho)s_1(x')\theta+s_2(x')-2\rho p(x')]
\right\}
	\end{gather*}
where
\begin{gather*}s_1(x'):=\sum_i(x'_{i1}+x'_{i2}),\quad s_2(x'):=\sum_i({x'_{i1}}^2+{x'_{i2}}^2),\quad p(x')=\sum_ix'_{i1}x'_{i2},\\
K_1(x')=\left[\int_{\mathbb R}f_{n,x'}(\theta)dQ(\theta)\right]^{-1}.
\end{gather*}
Note that, writing $c_n(\rho,\sigma,\tau)=[\tau\sqrt{2\pi}(2\pi\sigma^2\sqrt{1-\rho^2})^n]^{-1}$,
\begin{gather*}
\int_{\mathbb R}f_{n,x'}(\theta)dQ(\theta)=\\
c_n(\rho,\sigma,\tau)\int_{\mathbb R}
\left\{-\frac1{2\sigma^2(1-\rho^2)}[2n(1-\rho)\theta^2-2(1-\rho)s_1(x')\theta+s_2(x')-2\rho p(x')]-\frac1{2\tau^2}(\theta-\mu)^2]\right\}d\theta\\
=
c_n(\rho,\sigma,\tau)\int_{\mathbb R}
\exp\{-(A_1\theta^2-B_1(x')\theta+C_1(x'))\}d\theta\\
=\frac{\exp\left\{-\left(C_1(x')-\frac{B_1^2(x')}{4A_1}\right)\right\}}{\tau\sqrt{2\pi}[2\pi\sigma^2\sqrt{1-\rho^2}]^n}\int_{\mathbb R}\exp\left\{-A_1\left(\theta-\frac{B_1(x')}{2A_1}\right)^2\right\}d\theta\\
=\frac{\exp\left\{-\left(C_1(x')-\frac{B_1^2(x')}{4A_1}\right)\right\}}{\tau\sqrt{2A_1}[2\pi\sigma^2\sqrt{1-\rho^2}]^n}
	\end{gather*}
where
\begin{gather*}
A_1=\frac{n}{\sigma^2(1+\rho)}+\frac1{2\tau^2},\quad
B_1(x')=\frac{s_1(x')}{\sigma^2(1+\rho)}+\frac{\mu}{\tau^2},\quad C_1(x')=\frac{s_2(x')-2\rho p(x')}{2{\sigma^2(1-\rho^2)}}+\frac{\mu^2}{2\tau^2}.
	\end{gather*}
Hence
$$K_1(x')=\tau\sqrt{2A_1}[2\pi\sigma^2\sqrt{1-\rho^2}]^n\exp\left\{C_1(x')-\frac{B_1^2(x')}{4A_1}\right\}
$$
The posterior predictive density given $x'\in(\mathbb R^2)^n$ is
\begin{gather*}
f_{n,x'}^*(x):=\int_{\mathbb R}f_\theta(x)r_{n,x'}^*(\theta)g(\theta)d\theta=
\frac{K_1(x')}{\tau\sqrt{2\pi}[2\pi\sigma^2\sqrt{1-\rho^2}]^{n+1}}\cdot\\
\int_{\mathbb R}
\exp\left\{
-\frac{[2(n+1)(1-\rho)\theta^2-2(1-\rho)(x_1+x_2+s_1(x'))\theta+(x_1^2+x_2^2+s_2(x'))-2\rho (x_1x_2+p(x'))]}{2\sigma^2(1-\rho^2)}\right.\\
\left.-\frac{(\theta^2-2\mu\theta+\mu^2)}{2\tau^2}
\right\}d\theta=\\
K_2(x')\int_{\mathbb R}\exp\{-(A_2\theta^2-B_2(x,x')\theta+C_2(x,x'))\}d\theta=\\
K_2(x')\exp\left\{-\left(C_2(x,x')-\frac{B_2^2(x,x')}{4A_2}\right)\right\}
\int_{\mathbb R}\exp\left\{-A_2\left(\theta-\frac{B_2(x,x')}{2A_2}\right)^2\right\}d\theta=\\
K_2(x')\sqrt{\frac{\pi}{A_2}}\exp\left\{-\left(C_2(x,x')-\frac{B_2^2(x,x')}{4A_2}\right)\right\}
	\end{gather*}	
where
\begin{gather*}	
	K_2(x')=\frac{K_1(x')}{\tau\sqrt{2\pi}[2\pi\sigma^2\sqrt{1-\rho^2}]^{n+1}}=\frac{\sqrt{A_1}\exp\left\{C_1(x')-\frac{B_1^2(x')}{4A_1}\right\}}{2\pi^{3/2}\sigma^2\sqrt{1-\rho^2}},\quad A_2=\frac{(n+1)(1+\rho)}{\sigma^2}+\frac1{2\tau^2}\\
	B_2(x.x')=\frac{x_1+x_2+s_1(x')}{\sigma^2(1+\rho)}+\frac{\mu}{\tau^2},\quad 
	C_2(x,x')=\frac{x_1^2+x_2^2+s_2(x')-2\rho(x_1x_2+p(x'))}{2\sigma^2(1-\rho^2)}+\frac{\mu^2}{\tau^2}
	\end{gather*}
We can write
\begin{gather*}	
C_2(x,x')-\frac{B_2^2(x,x')}{4A_2}=A_3x_1^2+A_3x_2^2+B_3x_1x_2+C_3(x')(x_1+x_2)+D_3(x')
	\end{gather*}	
where
\begin{gather*}	
A_3=\frac1{2\sigma^2(1-\rho^2)}-\frac{\tau^2}{\sigma^2(1+\rho)^2[4(n+1)(1+\rho)\tau^2+2\sigma^2]},\\
B_3=-\frac{\rho}{\sigma^2(1-\rho^2)}-\frac{\tau^2}{\sigma^2(1+\rho)^2[2(n+1)(1+\rho)\tau^2+\sigma^2]},\\
C_3(x')=-\frac{\tau^2s_1(x')+\mu\sigma^2(1+\rho)}{\sigma^2(1+\rho)^2[2(n+1)(1+\rho)\tau^2+\sigma^2]},\\
D_3(x')=\frac{\tau^2s_2(x')-2\rho\tau^2p(x')+2\mu^2\sigma^2(1-\rho^2)}{2\tau^2\sigma^2(1-\rho^2)}-\frac{\tau^4s_1(x')^2+\mu^2\sigma^4(1+\rho)^2+2\mu\sigma^2(1+\rho)s_1(x')}{\tau^2\sigma^2(1+\rho)^2[4(n+1)(1+\rho)\tau^2+2\sigma^2]}
\end{gather*}
It is readily shown that $A_3>0$. It follows that the posterior predictive density $f^*_{n,x'}$ is the density of a normal bivariate distribution $N_2(m(x'),\Sigma)$ where
$$m(x')=\binom{m_1(x')}{m_2(x')},\quad \Sigma=
\sigma_1^2\left(\begin{array}{cc}
1 & \rho_1\\
\rho_1 & 1
\end{array}\right)
$$
being 
$$\rho_1=\frac{B_3}{2A_3},\quad \sigma_1^2=\frac{2A_3}{4A_3^2-B_3^2},\quad
m_1(x')=m_2(x')=-\frac{C_3(x')}{2(1-\rho_1)}.
$$
It is easy to see that, as can be expected, $|\rho_1|<1$ and $\sigma_1^2>0$. 
If we denote
$$a_n(\rho,\sigma,\tau):=2(n+1)(1+\rho)+\frac{\sigma^2}{\tau^2}
$$
We can write
\begin{gather*}	
\rho_1=-\frac{a_n(\rho,\sigma,\tau)+\frac{1-\rho}{1+\rho}}{a_n(\rho,\sigma,\tau)-\frac{1-\rho}{1+\rho}}\cdot\rho,\quad 
\sigma_1^2=\frac{a_n(\rho,\sigma,\tau)}{a_n(\rho,\sigma,\tau)-\frac{1-\rho}{1+\rho}}\cdot\sigma^2,\\
m_1(x')=m_2(x')=\frac{s_1(x')+(1+\rho)\frac{\sigma^2}{\tau^2}\mu}{2(1-\rho_1)(1+\rho)^2\sigma^2a_n(\rho,\sigma,\tau)}
\end{gather*}

It follows that the conditional distribution 
$$\left({R^*_{n,x'}}^{\!\!\!\!\!R}\right)^{\pi_2|\pi_1=x_1} \!\!:=N\big((1-\rho_1)m_1(x')+\rho_1x_1,\sigma_1^2(1-\rho_1^2)\big)
$$
is the Bayes estimator of the conditional distribution 
$$P_\theta^{X_2|X_1=x_1}=N\big((1-\rho)\theta+\rho x_1,\sigma^2(1-\rho^2)\big)
$$
for the squared total variation function, 
 and its density ${f^*_{n,x'}}^{\!\!\!\!\!\pi_2|\pi_1=x_1}$ is the Bayes estimator of the conditional density
$$f_\theta^{X_2|X_1=x_1}(x_2)=\frac1{\sigma\sqrt{2\pi(1-\rho^2)}}
\exp\left\{-\frac{1}{2\sigma^2(1-\rho^2)}[x_2-(1-\rho)\theta-\rho x_1]^2\right\}
$$ 
for the $L^1$-squared loss function. 

Moreover, its mean
$$E_{{R^*_{n,x'}}^{\!\!\!\!\!R}}(\pi_2|\pi_1=x_1)=(1-\rho_1)m_1(x')+\rho_1x_1
$$
is the Bayes estimator of the regression curve 
$$E_\theta(X_2|X_1=x_1)=(1-\rho)\theta+\rho x_1
$$
for the squared error loss function. 
	\end{myexa}

\section{Acknowledgements.}
This paper has been supported by the Junta de Extremadura (Spain) under the grant Gr18016.
\vspace{1ex}

\section{References.}

\begin{itemize}

	\item Barra, J.R. (1971) Notions Fondamentales de Statistique Math\'ematique, Dunod, Paris.

	\item Bogachev, V.I. (2007) Measure Theory, Vol. II, Springer, Berlin.

	\item Florens, J.P., Mouchart, M.,  Rolin, J.M. (1990) Elements of Bayesian Statistics, Marcel Dekker, New York.
	
	\item Geisser, S. (1993) Predictive Inference: An Introduction, Springer Science+ Business Media, Dordrecht.
	
	\item Gelman, A., Carlin, J.B., Stern, H.S., Dunson, D.B., Vehtari, A., Rubin, D.B. (2014) Bayesian Data Analysis, 3rd ed., CRC Press.

	\item Ghosal, S., Vaart, A.v.d. (2017) Fundamentals of Noparametric Bayesian Inference, Cambridge University Press, Cambridge UK. 

	\item Lehmann, E.L., Casella, G. (1998) Theory of Point Estimation, Second Edition, Springer, New York.

    \item Nadaraya, E.A. (1989) Nonparametric Estimation of Probability Densities and Regression Curves, Kluwer Academic Publishers, Dordrecht. 
    
    \item Nogales, A.G. (2021) On Bayesian Estimation of Densities and Sampling Distributions: the Posterior Predictive Distribution as the Bayes Estimator, to appear in Statistica Neerlandica.

	\item Rubin, D.B. (1984) Bayesianly justifiable and relevant frequency calculations for the applied statisticians, The Annals of Statistics, 12(4), 1151-1172.

\end{itemize}

\end{document}